\documentclass{proc-l}

\usepackage{amsfonts}
\usepackage{amsmath}
\usepackage{amssymb}

\newtheorem{theorem}{Theorem}[section]
\newtheorem{lemma}[theorem]{Lemma}

\theoremstyle{definition}

\theoremstyle{remark}

\numberwithin{equation}{section}

\begin{document}
\title[Global solutions]{Global solutions to special Lagrangian equations}
\author{Yu YUAN}
\address{Department of Mathematics, Box 354350\\
University of Washington\\
Seattle, WA 98195}
\email{yuan@math.washington.edu}
\thanks{The author is partially supported by an NSF grant and a Sloan Research
Fellowship. The author was a visiting fellow at the Australian National
University while this work was done.}

\commby{Jon Wolfson}

\begin{abstract}
We show that any global solution to the special Lagrangian equations with the
phase larger than a critical value must be quadratic.

\end{abstract}
\maketitle

\section{Introduction}

In this note, we show that any global solution $u$ in $\mathbb{R}^{n}$ to the
special Lagrangian equation%

\begin{equation}
\sum_{i=1}^{n}\arctan\lambda_{i}=c \label{E}%
\end{equation}
with phase $\left\vert c\right\vert >\frac{\pi}{2}\left(  n-2\right)  $ must
be a quadratic polynomial, where $\lambda_{i}s$ are the eigenvalues of the
Hessian $D^{2}u.$ Recall the Bernstein theorem, any global solution to the
minimal surface equation $div\left(  \frac{\nabla u}{\sqrt{1+\left\vert \nabla
u\right\vert ^{2}}}\right)  =0$ in $\mathbb{R}^{7}$ must be a linear function.

Equation (\ref{E}) stems from the special Lagrangian geometry [HL]. The
Lagrangian graph $\left(  x,\nabla u\left(  x\right)  \right)  \subset
\mathbb{R}^{n}\times\mathbb{R}^{n}$ is called special when the phase or the
argument of the complex number $\left(  1+\sqrt{-1}\lambda_{1}\right)
\cdots\left(  1+\sqrt{-1}\lambda_{n}\right)  $ is constant $c,$ that is, $u$
satisfies equation (\ref{E}), and it is special if and only if $\left(
x,\nabla u\left(  x\right)  \right)  $ is a minimal surface in $\mathbb{R}%
^{n}\times\mathbb{R}^{n}$ [HL, Theorem 2.3, Proposition 2.17]. To be precise,
we state

\begin{theorem}
Let $u$ be a smooth solution in $\mathbb{R}^{n}$ to (\ref{E} ) with
$\left\vert c\right\vert >\frac{\pi}{2}\left(  n-2\right)  ,$ then $u$ is quadratic.
\end{theorem}

Fu [F] proved Theorem 1.1 in the two dimensional case. Indeed (\ref{E}) with
$c=\frac{\pi}{2}$ in 2-d case also has the Monge-Amp\`{e}re form $\det
D^{2}u=1$ and J\"{o}rgens already showed Theorem 1.1 in this special case
earlier on (cf. [N]).

Other Bernstein/Liouville type results concerning (\ref{E}) are in order.
Borishenko [B] showed that any convex solution with linear growth to (\ref{E})
with $c=k\pi$ is linear. The author [Y] proved that any convex solution to
(\ref{E}) must be quadratic. For $c=k\pi$ in $n=3$ case, (\ref{E}) has another
form%
\begin{equation}
\bigtriangleup u=\det D^{2}u. \label{El-m}%
\end{equation}
It was proved in [BCGJ] that any strictly convex solution to (\ref{El-m}) with
quadratic growth must be quadratic. Under the assumption that the Hessian is
bounded and $\lambda_{i}\lambda_{j}\geq-\frac{3}{2},$ it was also showed in
[TW] that any global solution to (\ref{E}) is quadratic.

The heuristic idea of the proof of Theorem 1.1 is to find a better graph
representation of $\left(  x,\nabla u\left(  x\right)  \right)  $ so that the
Hessian of the new potential is bounded and the new potential function
satisfies a convex uniformly elliptic equation. By Krylov-Evan's [K] [E]
interior H\"{o}lder estimates on the Hessian, we draw the conclusion.

As there are nontrivial global harmonic solutions to (1.1) with $c=0$ in case
$n=2,$ we guess (1.1) with $c=\frac{\pi}{2}\left(  n-2\right)  $ also has
nontrivial global solutions in the higher dimensional case. Observe in the
case $n=3$ and $c=\frac{\pi}{2}$, (\ref{E}) also takes the interesting form
$\lambda_{1}\lambda_{2}+\lambda_{2}\lambda_{3}+\lambda_{3}\lambda_{1}=1.$

\section{Proof}

Step1. We first find a better graph representation of $M$ through Lewy
rotation (cf. [N]) so that the Hessian of the potential function is bounded.
By symmetry we only consider the case $c>\frac{\pi}{2}\left(  n-2\right)  .$
Let $\sum_{i=1}^{n}\theta_{i}=\frac{\pi}{2}\left(  n-2\right)  +\delta,$ where
$\theta_{i}=\arctan\lambda_{i}\in\left(  -\frac{\pi}{2},\frac{\pi}{2}\right)
$ and $\delta\in\left(  0,\pi\right)  .$ Note that%
\begin{equation}
-\frac{\pi}{2}+\frac{\left(  n-1\right)  }{n}\delta<\theta_{i}-\frac{\delta
}{n}<\frac{\pi}{2}-\frac{\delta}{n}.\label{Ehessianbound}%
\end{equation}
The first inequality follows from $\frac{\pi}{2}\left(  n-2\right)
+\delta<\theta_{i}+\frac{\pi}{2}\left(  n-1\right)  .$ We rotate the $\left(
x,y\right)  \in\mathbb{R}^{n}\times\mathbb{R}^{n}$ coordinate system to
$\left(  \bar{x},\bar{y}\right)  $ by $\frac{\delta}{n},$ namely, $\bar
{x}=x\cos\frac{\delta}{n}+y\sin\frac{\delta}{n},$ $\bar{y}=-x\sin\frac{\delta
}{n}+y\cos\frac{\delta}{n}.$ In terms of complex variables $z=x+\sqrt{-1}y$,
that is, we identify $\mathbb{R}^{n}\times\mathbb{R}^{n}\ $with $\mathbb{C}%
^{n},$ the rotation takes the form $\bar{z}=e^{-\sqrt{-1}\delta/n}z.$ Then $M$
has a new parametrization%
\[%
\genfrac{\{}{.}{0pt}{}{\bar{x}=x\cos\frac{\delta}{n}+\nabla u\left(  x\right)
\sin\frac{\delta}{n}}{\bar{y}=-x\sin\frac{\delta}{n}+\nabla u\left(  x\right)
\cos\frac{\delta}{n}}%
.
\]
By (\ref{Ehessianbound}), $M=\left(  x,\nabla u\left(  x\right)  \right)  $ is
still a graph over the whole $\bar{x}$ space $\mathbb{R}^{n}.$ Further the
rotation belongs to $U\left(  n\right)  ,$ then $M$ is also a special
Lagrangian graph $\left(  \bar{x},\nabla\bar{u}\left(  \bar{x}\right)
\right)  ,$ where $\bar{u}$ is a smooth function [HL, p.87, Proposition 2.17].
Let $\bar{\lambda}_{i}$ be the eigenvalues of $D^{2}\bar{u},$ then
$\bar{\theta}_{i}=\arctan\bar{\lambda}_{i}=\theta_{i}-\frac{\delta}{n}%
\in\left(  -\frac{\pi}{2}+\frac{\left(  n-1\right)  }{n}\delta,\frac{\pi}%
{2}-\frac{\delta}{n}\right)  .$ That is%
\[
\left\vert D^{2}\bar{u}\right\vert \leq C\left(  \delta\right)  .
\]
Finally $\bar{u}$ satisfies the equation%
\begin{equation}
\sum_{i=1}^{n}\arctan\bar{\lambda}_{i}=\frac{\pi}{2}\left(  n-2\right)
.\label{Erotation}%
\end{equation}

Step 2. We proceed with the following lemma, which is Lemma 8.1 in [CNS] when
$n$ is even and $c=\frac{\pi}{2}\left(  n-2\right)  .$

\begin{lemma}
Let $f\left(  \lambda_{1},\cdots,\lambda_{n}\right)  =\sum_{i=1}^{n}%
\arctan\lambda_{i}$ and $\Gamma=\left\{  \lambda|f\left(  \lambda\right)
=c\right\}  $ with $\left\vert c\right\vert \geq\frac{\pi}{2}\left(
n-2\right)  ,$ then $\Gamma$ is convex.
\end{lemma}

\begin{proof}
\label{of Lemma}We skip the case $n=1.$ By symmetry we just consider the case
$c\geq0.$ Set $c=\frac{\pi}{2}\left(  n-2\right)  +\delta$ with $\delta
\in\lbrack0,\pi).$ We may assume that $\theta_{i}=\arctan\lambda_{i}\geq0$ for
$i=1,\cdots,n-1.$ The normal of $\Gamma$ is $\nabla f=\left(  \cos^{2}%
\theta_{1},\cdots,\cos^{2}\theta_{n}\right)  .$ Let%
\[
A\triangleq-\frac{1}{2}D^{2}f=\left[
\begin{array}
[c]{ccc}%
\tan\theta_{1}\cos^{4}\theta_{1} &  & \\
& \cdots & \\
&  & \tan\theta_{n}\cos^{4}\theta_{n}%
\end{array}
\right]  .
\]
Take any tangent vector $T=\left(  t_{1},\cdots,t_{n}\right)  \in T_{\lambda
}\Gamma$, that is
\[
\sum_{i=1}^{n}t_{i}\cos^{2}\theta_{i}=0.
\]
We show that $A\left(  T,T\right)  \geq0.$

Case a) $\theta_{n}\geq0.$ Certainly it is true.

Case b) $\theta_{n}<0.$ First we know that $\theta_{i}>0$ for $i=1,\cdots,n-1$
and $\delta<\frac{\pi}{2}.$ Next we have%
\[
A\left(  T,T\right)  =\sum_{i=1}^{n-1}\tan\theta_{i}\cos^{4}\theta_{i}%
t_{i}^{2}+\tan\theta_{n}\cos^{4}\theta_{n}t_{n}^{2}.
\]
Now we use the trick in [CNS, p.299],%
\begin{align*}
\left(  -t_{n}\cos\theta_{n}\right)  ^{2}  &  =\left(  \sum_{i=1}^{n-1}%
t_{i}\cos^{2}\theta_{i}\right)  ^{2}\\
&  \leq\left(  \sum_{i=1}^{n-1}t_{i}^{2}\cos^{4}\theta_{i}\tan\theta
_{i}\right)  \left(  \sum_{i=1}^{n-1}\cot\theta_{i}\right)
\end{align*}
then%
\[
\tan\theta_{n}\cos^{4}\theta_{n}t_{n}^{2}\geq\left(  \sum_{i=1}^{n-1}t_{i}%
^{2}\cos^{4}\theta_{i}\tan\theta_{i}\right)  \left(  \sum_{i=1}^{n-1}%
\cot\theta_{i}\right)  \tan\theta_{n},
\]
and%
\begin{align*}
A\left(  T,T\right)   &  \geq\left(  \sum_{i=1}^{n-1}t_{i}^{2}\cos^{4}%
\theta_{i}\tan\theta_{i}\right)  \left[  1+\left(  \sum_{i=1}^{n-1}\cot
\theta_{i}\right)  \tan\theta_{n}\right] \\
&  =\left(  \sum_{i=1}^{n-1}t_{i}^{2}\cos^{4}\theta_{i}\tan\theta_{i}\right)
\left(  \sum_{i=1}^{n}\cot\theta_{i}\right)  \tan\theta_{n}.
\end{align*}
Let $\alpha_{i}=\frac{\pi}{2}-\theta_{i},$ we have%
\begin{align*}
\frac{\pi}{2}  &  <\pi-\delta=\alpha_{1}+\cdots+\alpha_{n}<\pi,\\
0  &  <\alpha_{1},\cdots,\alpha_{n-1}<\frac{\pi}{2}<\alpha_{n},
\end{align*}
and%
\[
\sum_{i=1}^{n}\cot\theta_{i}=\sum_{i=1}^{n-1}\tan\alpha_{i}+\tan\alpha_{n}.
\]
It follows that $\tan\alpha_{n}<0$ and%
\[
\frac{\tan\left(  \alpha_{1}+\cdots+\alpha_{n-1}\right)  +\tan\alpha_{n}%
}{1-\tan\left(  \alpha_{1}+\cdots+\alpha_{n-1}\right)  \tan\alpha_{n}}%
=\tan\left(  \alpha_{1}+\cdots+\alpha_{n}\right)  <0.
\]
Then $\tan\left(  \alpha_{1}+\cdots+\alpha_{n-1}\right)  +\tan\alpha_{n}<0.$
Note that $\alpha_{1}+\cdots+\alpha_{n-1}=\pi-\delta-\alpha_{n}<\frac{\pi}%
{2},$ we have%
\begin{align*}
\tan\left(  \alpha_{1}+\cdots+\alpha_{n-1}\right)   &  \geq\tan\alpha_{1}%
+\tan\left(  \alpha_{2}+\cdots+\alpha_{n-1}\right) \\
&  \cdots\\
&  \geq\tan\alpha_{1}+\tan\alpha_{2}+\cdots+\tan\alpha_{n-1}.
\end{align*}
So $\sum_{i=1}^{n}\cot\theta_{i}=\sum_{i=1}^{n}\tan\alpha_{i}$ $<0$ and
$A\left(  T,T\right)  \geq0.$ Therefore $\Gamma$ is convex (w.r.t. the normal
$\nabla f$ ).
\end{proof}

\textbf{Remark. }The level set $\Gamma$ is no longer convex nor concave when
$\left\vert c\right\vert <\frac{\pi}{2}\left(  n-2\right)  .$

\bigskip

Step 3. The final argument is standard. We now have global solution $\bar{u}$
with bounded Hessian on the convex level set $\Gamma,$ more precisely a convex
level set in the symmetric matrix space (cf. [CNS, p.276). In another word,
$\bar{u}$ satisfies (\ref{Erotation}), which is uniformly elliptic now. By
Krylov-Evans theorem ([K],[E])%
\[
\left[  D^{2}\bar{u}\right]  _{C^{\beta}\left(  B_{r}\right)  }\leq C\left(
n,\delta\right)  \frac{\left\Vert D^{2}\bar{u}\right\Vert _{L^{\infty}\left(
B_{2r}\right)  }}{r^{\beta}}\leq\frac{C\left(  n,\delta\right)  }{r^{\beta}},
\]
where $\beta=\beta\left(  n,\delta\right)  \in\left(  0,1\right)  .$ Let $r$
go to $+\infty,$ we see that $D^{2}\bar{u}$ is a constant matrix. Thus
$\left(  \bar{x},\nabla\bar{u}\right)  $ is a plane and consequently $u$ is quadratic.

\end{document}